\documentclass[a4paper,11pt,reqno]{amsart}
\usepackage{a4wide}
\usepackage{amsfonts,amsmath,amssymb}
\usepackage{graphicx}
\usepackage{url}
\usepackage{cite}

\newtheorem{thm}{Theorem}

\newtheorem{lem}[thm]{Lemma}

\theoremstyle{remark}
\newtheorem{rem}{Remark}

\begin{document}

\thispagestyle{plain}

\title{The number of fixed points of Wilf's partition involution}

\author[Stephan Wagner]{Stephan Wagner}
\address{Stephan Wagner \\
Department of Mathematical Sciences \\
Stellenbosch University \\
Private Bag X1 \\
Matieland 7602 \\
South Africa
}
\email{swagner@sun.ac.za}
\thanks{This material is based upon work supported by the National Research Foundation of South Africa under grant number 70560.}

\date{\today}


\begin{abstract}
Wilf partitions are partitions of an integer $n$ in which all nonzero multiplicities are distinct. On his webpage, the late Herbert Wilf posed the problem to find ``any interesting theorems'' about the number $f(n)$ of those partitions. Recently, Fill, Janson and Ward (and independently Kane and Rhoades) determined an asymptotic formula for $\log f(n)$. Since the original motivation for studying Wilf partitions was the fact that the operation that interchanges part sizes and multiplicities is an involution on the set of Wilf partitions, they mentioned as an open problem to determine a similar asymptotic formula for the number of fixed points of this involution, which we denote by $F(n)$. In this short note, we show that the method of Fill, Janson and Ward also applies to $F(n)$. Specifically, we obtain the asymptotic formula $\log F(n) \sim \frac12 \log f(n)$. 
\end{abstract}
\maketitle

\section{Introduction}\label{intro}

By a \emph{Wilf partition}, we mean an integer partition in which all nonzero multiplicities are distinct. On his webpage, the late Herbert Wilf posted a set of unsolved problems, among them the task to find ``any interesting theorems'' about the number $f(n)$ of such partitions. An important contribution to this problem was made recently, when Fill, Janson and Ward \cite{fill2012partitions} (and independently Kane and Rhoades \cite{kane2013asymptotics}) proved an asymptotic formula for $\log f(n)$:
\begin{equation}\label{eq:fjw}
\log f(n) \sim \frac13 (6n)^{1/3} \log n.
\end{equation}
Wilf's motivation for considering partitions with distinct multiplicities was the fact that the operation that interchanges parts and multiplicities is an involution on the set of all such partitions. For example, if we consider the Wilf partition
$$5 + 5 + 4 + 3 + 3 + 3 + 1+1+1+1+1+1$$
of $29$, then by interchanging parts and multiplicities, we get a new Wilf partition of $29$:
$$6 + 3 + 3 + 3 + 2 + 2 + 2 + 2 + 2 + 1 + 1 + 1 + 1.$$
Some Wilf partitions are fixed points of this involution, such as
\begin{equation}\label{wilfexample}
6 + 4 + 4 + 3 + 3 + 3 + 2 + 2 + 2 + 2 + 1 + 1 + 1 + 1 + 1 + 1,
\end{equation}
which is mapped to itself when parts and multiplicities are interchanged. We are interested in the number $F(n)$ of such fixed points, in particular asymptotics for it, which was left as an open problem by Fill, Janson and Ward. The first few terms of the sequences $f(n)$ (that counts all Wilf partitions) and $F(n)$ (Sloane's A098859 and A217605 respectively \cite{sloane2013encyclopedia}) are given in the following table:

\begin{table}[htbp]
\begin{center}
\begin{tabular}{|c|cccccccccccccccccccc|}
\hline
$n$ & 1 & 2 & 3 & 4 & 5 & 6 & 7 & 8 & 9 & 10 & 11 & 12 & 13 & 14 & 15 & 16 & 17 & 18 & 19 & 20 \\
\hline
$f(n)$ & 1 & 1 & 2 & 2 & 4 & 5 & 7 & 10 & 13 & 15 & 21 & 28 & 31 & 45 & 55 & 62 & 82 & 105 & 116 & 153 \\
\hline
$F(n)$ & 1 & 0 & 0 & 2 & 1 & 1 & 0 & 1 & 1 & 3 & 0 & 3 & 3 & 3 & 0 & 4 & 3 & 2 & 1 & 6 \\
\hline
\end{tabular}
\end{center}
\caption{Values of $F(n)$ for small $n$.}\label{table:fvalues}
\end{table}

The behaviour of $F(n)$ is quite erratic, much more so than for $f(n)$ (which is conjectured to be increasing). However, there is still an asymptotic formula for $\log F(n)$, as we show in the following:

\begin{thm}\label{thm:main}
The number $F(n)$ of fixed points of the involution on Wilf partitions satisfies
\begin{equation}\label{eq:main}
\log F(n) \sim \frac12 \log f(n) \sim \frac16 (6n)^{1/3} \log n.
\end{equation}
\end{thm}

The proof of this theorem follows very closely the ideas of Fill, Janson and Ward. It procedes in two stages: first we determine an upper bound for $F(n)$, then we construct a set of fixed points of the Wilf involution that is large enough to provide us with a matching lower bound. It is perhaps interesting to compare~\eqref{eq:main} to a similar combinatorial result: \emph{involutions} are permutations with only $1$- and $2$-cycles, much like Wilf partitions, where the part-multiplicity pairs form $1$- and $2$-cycles. If $I(n)$ is the number of involutions, then
$$\log I(n) \sim \frac12 \log n! \sim \frac{n}{2} \log n,$$
see \cite[Example VIII.5]{flajolet2009analytic}. Thus when $I(n)$ is compared to the total number $n!$ of permutations, one encounters a similar factor of $\frac12$.

\section{Proof of the main theorem}

\subsection{The upper bound}

We represent partitions as sequences of part-multiplicity pairs $(p_i,m_i)$ -- each such pair represents $m_i$ terms in the partition that are equal to $p_i$. If a Wilf partition is a fixed point of our involution and the part-multiplicity pair $(p_i,m_i)$ occurs in it, then so does the pair $(m_i,p_i)$. 

We split such a Wilf partition into two parts:
\begin{itemize}
\item The part that consists of all part sizes that are equal to their multiplicities, i.e., all part-multiplicity pairs $(p_i,m_i)$ such that $p_i = m_i$.
\item The rest, which consists of pairs of part-multiplicity pairs of the form $(p_i,m_i)$ and $(m_i,p_i)$ with $m_i \neq p_i$.
\end{itemize}
For instance, the partition in~\eqref{wilfexample} is split into the part-multiplicity pair $(3,3)$ and the rest, which consists of the pairs $(6,1)$, $(4,2)$, $(2,4)$ and $(1,6)$. The first part is equivalent to a partition of some number $k \leq n$ into distinct squares. Let $q(k)$ denote the number of partitions of $k$ into squares. It is known (as a special case of the Meinardus scheme, see \cite[Chapter 6]{andrews1998theory}) that
$$\log q(k) \sim C k^{1/3}$$
for a constant $C = 3 (\sqrt{\pi} \zeta(3/2)/4)^{2/3}$. If we take the remaining part-multiplicity pairs $(p_i,m_i)$ and remove all those for which $p_i < m_i$, then we obtain a Wilf partition of $(n-k)/2$. Moreover, as it was shown by Fill, Janson and Ward, a Wilf partition of $n$ can have at most $(6n)^{1/3}$ distinct parts: if there are $r$ distinct parts, then the multiplicities are at least $1,2,\ldots,r$, and so we have
$$n \geq 1 \cdot r + 2 \cdot (r-1) + \ldots + r \cdot 1 = \frac{r(r+1)(r+2)}{6} \geq \frac{r^3}{6}.$$
Since we are cutting in half, the remaining Wilf partition has at most $\frac12 (6(n-k))^{1/3}$ distinct parts.  A Wilf partition of $n$ with $r$ distinct parts can be obtained by the following two-step process (cf. \cite{fill2012partitions}):
\begin{itemize}
\item start with a partition of $n$ into $r$ parts,
\item split every part $x_i$ into $m_i$ copies of $p_i$ such that $p_i m_i = x_i$. If we require in addition that $p_i > m_i$, then there are at most $d(x_i)/2$ possibilities for the choice of $p_i$ and $m_i$, where $d(\cdot)$ denotes the number of divisors.
\end{itemize}
Note that not every partition that we get from this process is actually a Wilf partition, but every Wilf partition can be obtained in this way. As in \cite{fill2012partitions}, we 
now make use of the estimate
\begin{equation}\label{eq:partestim}
p(m,r) = O \left( \frac{m^{r-1}}{r!(r-1)!} \right)
\end{equation}
for the number of partitions of $m$ into $r$ summands, which holds for $r = O(n^{1/3})$ (cf. \cite[Section 7.2.1.4, Exercise 34]{knuth2005art}): $p(m,r)$ is also the number of partitions of $m + r(r-1)/2$ into \emph{distinct} parts (add $0,1,2,\ldots,r-1$ to the parts to make them distinct). Each such partition can be turned into a \emph{composition} by permuting the $r$ parts in all possible $r!$ ways. It is well known \cite[Example I.6]{flajolet2009analytic} that the number of compositions of $s$ into $r$ parts is $\binom{s-1}{r-1}$, thus
$$p(m,r) \leq \frac{1}{r!} \binom{m+r(r-1)/2-1}{r-1} \leq \frac{(m+r^2/2)^{r-1}}{r!(r-1)!},$$
from which~\eqref{eq:partestim} follows. Moreover, we need the estimate
$$\log d(n) = O \left( \frac{\log n}{\log \log n} \right)$$
for the number of divisors of $n$ (see \cite[Theorem 317]{hardy1979introduction}).
These give us the following upper bound:
\begin{align*}
F(n) \leq \sum_{\substack{k=0 \\ k \equiv n \bmod 2}}^n \sum_{r \leq (6(n-k))^{1/3}/2} q(k) p((n-k)/2,r) \max_{ x \leq n } (d(x)/2)^r \\
\leq n^2 q(n) \max_{m \leq n/2, r \leq (6n)^{1/3}/2} p(m,r) \max_{ x \leq n } (d(x)/2)^r
\end{align*}
and thus
$$\log F(n) \leq \frac16 (6n)^{1/3} \log n \left( 1 + O \left( \frac{1}{\log \log n} \right) \right) = f(n) (1+o(1))/2.$$

\subsection{The lower bound}

It remains to prove an inequality in the opposite direction. To this end, we employ the ideas of Fill, Janson and Ward once again. Our goal is to construct a sufficiently rich set of fixed points of the involution. This is slightly more difficult than the first part of the proof and requires the following lemma:

\begin{lem}\label{lem:sufficient_set}
There exists a set $\mathcal{A} = \{a_1,a_2,a_3,\ldots\}$ of positive integers such that
\begin{itemize}
\item Every positive integer, except for $2,3,7,11$ and $15$, has a Wilf partition that is a fixed point of the involution and only uses part sizes and multiplicities in $\mathcal{A}$.
\item $|\mathcal{A} \cap [1,m]| = O(\log m)$.
\end{itemize}
\end{lem}

\begin{rem}
In particular, this means that $2,3,7,11$ and $15$ are the only positive integers that do not have a Wilf partition that is a fixed point (it is easy to check that there is indeed none in those five cases -- see also Table~\ref{table:fvalues}).
\end{rem}

\begin{proof}
We set $a_k = k$ for $1 \leq k \leq 10$. One can verify directly that all positive integers $\leq 136$, except for $2,3,7,11$ and $15$, have Wilf partitions that are fixed points of the involution and only use $1,2,\ldots,10$ as part sizes and multiplicities. Now set $b_{10} = 136$ and define $a_k$ and $b_k$ for $k > 10$ recursively as follows:
$$a_k = \lfloor \sqrt{b_{k-1}-15} \rfloor \quad \text{and} \quad b_k = b_{k-1} + a_k^2.$$
It is easy to see that $a_k = \Omega((\sqrt{2}-\epsilon)^k)$ for any $\epsilon > 0$, so the second statement holds. For the first statement, we show that $a_1,\ldots,a_k$ are sufficient to obtain a suitable Wilf partition of any integer in the interval $[16,b_k]$. This is true for $k=10$ by our choice of $b_{10}$. For $k > 10$, the induction hypothesis guarantees the existence of a suitable Wilf partition for all integers in the interval $[16,b_{k-1}]$. By adding a part-multiplicity pair $(a_k,a_k)$ if necessary, we can also cover all values from $16+a_k^2 \leq b_{k-1}+1$ to $b_{k-1} + a_k^2 = b_k$, which completes the induction.
\end{proof}

Now we complete the proof of our main theorem by providing a lower bound for $F(n)$. Fix a large integer $K$, and let $C$ be a constant such that $|\mathcal{A} \cap [1,m]| \leq C \log m$ for sufficiently large $m$, as is guaranteed by Lemma~\ref{lem:sufficient_set}.  Now let $R$ be the integer nearest to
$$\frac{1}{2K} \left( \frac{6n}{1+3/K} \right)^{1/3}.$$
Consider the smallest $2RK$ integers that are not elements of $\mathcal{A}$, and denote them by $x_1,x_2,\ldots,$ $x_{2RK}$. By our choice of $R$, those integers are all smaller than $n$ for sufficiently large $n$, which means that there are at most $C \log n$ elements of $\mathcal{A}$ that are less than $x_{2RK}$. This implies that
$$x_j \leq C \log n + j.$$
Choose K permutations of $\{1,2,\ldots,R\}$ independently, and call them $\sigma_1,\ldots,\sigma_K$. Now we construct a Wilf partition of $n$ that is also a fixed point of the involution as follows:
\begin{itemize}
\item For $1 \leq i \leq K$ and $1 \leq j \leq R$, we include $x_{(i-1)R+j}$ copies of $x_{(2K-i)R+\sigma_i(j)}$, and vice versa.
\item The total contribution of these parts is
\begin{align*}
m &= 2 \sum_{i=1}^K \sum_{j=1}^R x_{(i-1)R+j} \cdot  x_{(2K-i)R+\sigma_i(j)} \\
&\leq 2 \sum_{i=1}^K \sum_{j=1}^R (C \log n + iR)(C \log n + (2K+1-i)R) \\
&= 2 C^2 KR \log^2 n + 2 CK(2K+1)R^2\log n + \frac{2K(K+1)(2K+1)}{3}R^3 \\
&= \frac{2K^2+3K+1}{2K^2+6K} n + O \left( n^{2/3} (K + \log n) \right)
\end{align*}
by our choice of $R$. For sufficiently large $n$, this is less than $n-16$. Thus we can find a Wilf permutation of $n - m$ using only parts and multiplicities in $\mathcal{A}$ that is also a fixed point of the involution. Since none of the $x_i$ lies in $\mathcal{A}$ by our choice, we can combine the two to a Wilf permutation of $n$ that still has the desired property.
\end{itemize}
It follows immediately that
$$F(n) \geq R!^K$$
and thus
$$\log F(n) \geq K \log R! = KR (\log R + O(1)) = \frac12 \left( \frac{6n}{1+3/K} \right)^{1/3} \left( \frac13 \log n + O(1) \right).$$
Since $K$ was arbitrary, we have
$$\log F(n) \geq \frac16 (6n)^{1/3} \log n (1+o(1)) = f(n)(1+o(1))/2,$$
which completes the proof of Theorem~\ref{thm:main}.

\begin{rem}
Setting e.g. $K = \lfloor \log n \rfloor$, we can get a slightly  more precise lower bound:
$$\log F(n) \geq \frac16 (6n)^{1/3} (\log n + O(\log \log n)).$$
\end{rem}

\bibliographystyle{abbrv}
\bibliography{wilfpartitions}

\end{document}